\documentclass[twoside,noinfoline,dvips,psamsfonts]{article}
\RequirePackage[lnms,amsthm,amsmath,natbib]{imsart}
\RequirePackage{graphicx}
\RequirePackage{hyperref}
\usepackage{amsthm,amsmath}
\usepackage{epsfig}
\usepackage{amssymb}
\usepackage{color}
\doi{10.1214/074921706000000022}
\pubyear{2006}
\volume{48}
\volumetitle{Dynamics \& Stochastics}
\firstpage{1}
\lastpage{15}

\startlocaldefs
\newtheorem{theorem}{Theorem}
\newtheorem{proposition}[theorem]{Proposition}
\newtheorem{corollary}[theorem]{Corollary}
\theoremstyle{definition}
\newtheorem{definition}[theorem]{Definition}
\newtheorem{example}{Example}
\newtheorem*{examp}{Example}
\newtheorem*{rem}{Remark}

\numberwithin{theorem}{section}
\numberwithin{equation}{section}
\endlocaldefs

\begin{document}

\begin{frontmatter}

\title{Polymer pinning in a random medium
      as influence percolation}
\runtitle{Polymer pinning in a random medium
      as influence percolation}

\begin{aug}
\author[A]{\fnms{V.} \snm{Beffara}\corref{}\ead[label=e1]{vbeffara@ens-lyon.fr}},
\author[B]{\fnms{V.} \snm{Sidoravicius}\ead[label=e2]{vladas@impa.br}},
\author[C]{\fnms{H.} \snm{Spohn}\ead[label=e3]{spohn@ma.tum.de}}
\and%
\author[D]{\fnms{M. E.} \snm{Vares}\ead[label=e4]{eulalia@cbpf.br}}

\affiliation{CNRS--ENS-Lyon, IMPA, TU-M\"{u}nchen and CBPF}
\address[A]{CNRS -- UMPA -- ENS Lyon,
46 All\'ee d'Italie,
F-69364 Lyon Cedex 07,
France,
\printead{e1}
}

\address[B]{IMPA,
Estrada Dona Castorina, 110,
Rio de Janeiro RJ 22460-320,
Brasil,
\printead{e2}
 }

\address[C]{Zentrum Mathematik,
TU Muenchen,
D-85478 Garching,
Germany,
\printead{e3}
 }

\address[D]{CBPF,
Rua Dr. Xavier Sigaud, 150,
Rio de Janeiro RJ 22290-180,
Brasil,
\printead{e4}
 }
\end{aug}

\contributor{Beffara, V.}{CNRS--ENS-Lyon}
\contributor{Sidoravicius, V.}{IMPA}
\contributor{Spohn, H.}{TU-M\"{u}nchen}
\contributor{Vares, M. E.}{CBPF}

\runauthor{V. Beffara {\it et al.}}

\begin{abstract}
  In this article we discuss a set of geometric ideas which shed some light on
  the question of  directed polymer pinning in the  presence of bulk disorder.
  Differing from standard methods and  techniques, we transform the problem to
  a particular dependent percolative  system and relate the pinning transition
  to a percolation transition.
\end{abstract}

\begin{keyword}[class=AMS]
\kwd[primary ]{60G55}
\kwd{60K35}
\kwd[; secondary ]{60G17}.
\end{keyword}

\begin{keyword}
\kwd{last-passage percolation}
\kwd{pinning}
\kwd{exclusion process}.
\end{keyword}

\end{frontmatter}

\section{Introduction}

\paragraph{Motivating example.}

The  totally asymmetric  exclusion process  (TASEP  for short)  is defined  as
follows on the set $\mathbb Z$: At time $0$, a (possibly random) configuration
of particles  is given,  in such  a way that  each site  contains at  most one
particle. To each  edge of the lattice is associated  a Poisson clock of
intensity $1$.   Whenever this  clock rings,  and there is  a particle  at the
left-end vertex  of this  edge and  no particle at  the right-end  vertex, the
particle moves to the right; otherwise  the ring of the clock is ignored.  The
product over  $\mathbb Z$ of  Bernoulli measures of density  $\rho\in(0,1)$ is
invariant  by this dynamics;  in that  case, the  average number  of particles
passing  through  the origin  up  to time  $t$  is  equal to  $\rho(1-\rho)t$,
\emph{i.e.} the flux through a given bond is exactly $\rho(1-\rho)$.

The process  is modified  at the  origin, by imposing  that the  Poisson clock
associated  with the  bond  $e_0=\langle 0,1  \rangle$  is $\lambda>0$.   When
$\lambda>1$, one can still prove that  the above expression for the flux holds
asymptotically, although  the Bernoulli measure is not  an equilibrium measure
anymore.  From now on, we shall assume that $0<\lambda\leq1$.

One of the  fundamental questions in driven flow is  to understand under which
conditions such  a static  obstruction results in  the formation of  a
``platoon''
starting at  the origin and propagating leftward.   A convenient quantitative
criterion  for platoon  formation  is to  start  the TASEP  with step  initial
conditions, i.e.~all sites $x\leq 0$ filled and all sites $x\geq 1$ empty, and
to  consider  the  average  current,  $j(\lambda)$, in  the  long  time  limit
$t\to\infty$.   $j(0)=0$,  $j$  is  non-decreasing, and  $j(\lambda)=1/4$  for
$\lambda=1$.   Thus   the  issue  is  to  determine   the  critical  intensity
$\lambda_\mathrm{c}$, which  is defined as  the supremum of all  the $\lambda$
for   which   $j(\lambda)<\frac{1}{4}$.     Estimates   for   the   value   of
$\lambda_\mathrm{c}$  were  given in  \cite{L}  and  \cite{CR},  and the  full
hydrodynamical picture is proved in \cite{S}.

The blockage problem  for the TASEP has been studied  numerically and by exact
enumerations.  On   the  basis  of   these  data,  in  \cite{JL2}   the  value
$\lambda_\mathrm{c}=1$  is   conjectured.   Recently  this   result  has  been
challenged (\cite{HTN}) and $\lambda_\mathrm{c}\cong 0.8$ is claimed.

One-dimensional  driven lattice  gases  belong to  the  universality class  of
Kardar-Parisi-Zhang (KPZ)  type growth  models.  In particular  the asymmetric
simple  exclusion process can  be represented  as the  so-called body-centered
solid-on-solid version  of (1+1)-dimensional  polynuclear growth model,  or as
directed polymer  subject to  a random  potential. The TASEP  also  has a well
known  representation  in  terms   of  last-passage  percolation  (or  maximal
increasing subsequences, known as Ulam's  problem).  In this article we choose
to work  in the setup of last  passage percolation.  The slow  bond induces an
extra   line    of   defects   relative   to   the    disordered   bulk.    If
$\lambda<\lambda_\mathrm{c}$,    the    optimal    path,   \emph{i.e.}     the
\emph{geodesic},    is    pinned    to    the    line    of    defects.     As
$\lambda\to\lambda_\mathrm{c}$, the geodesic  wanders further and further away
from the  line of defects  and the density  of intersections with the  line of
defects tends to zero.   For $\lambda>\lambda_\mathrm{c}$, the fluctuations of
the  geodesic  are  determined  by  the  bulk, and  the  line  of  defects  is
irrelevant.

In   the  present   work  we   will  not   establish  the   actual   value  of
$\lambda_\mathrm{c}$ (since we cannot) or settle the question as to whether it
is equal  to $1$;  the main goal  of this paper  is to  describe a new  way of
looking at the problem which gives some insight about the precise behavior of
the system.   More precisely, we show  how the problem can  be studied through
particular dependent  percolative systems constructed  in such a way  that the
pinning transition can be understood in terms of a percolation transition.

\section{Interpretation as pinning in Ulam's problem}

In this section we describe the representation of our initial problem in terms
of Ulam's problem or polynuclear  growth.  For more detailed explanations, see
\emph{e.g.} \cite{AD} or \cite{PS1}.

\subsection{The Model}

Let  $P^{(2)}$ be the  distribution of  a Poisson  point process  of intensity
$\lambda^{(2)}= 1$ in the plane  $\mathbb R^2$, and let $\Omega^{(2)}$ the set
of all its possible configurations; for  all $n>0$, let $P^{(2)}_n$ be the law
of  its  restriction to  the  square ${\cal  Q}_n  :=  [0,n]\times [0,n]$  and
$\Omega_n^{(2)}$ be the configuration space of the restricted process.

The \emph{maximal  increasing subsequence problem}, or Ulam's  problem, can be
formulated   in   the  following   geometric   way:   Given  a   configuration
$\omega^{(2)}\in\Omega^{(2)}$ and  its restriction $\omega_n^{(2)}$  to ${\cal
  Q}_n$, look  for an oriented path  $\pi$ (moving only  upward and rightward)
from $(0,0)$  to $(n,n)$  collecting as many  points from  $\omega_n^{(2)}$ as
possible;  as  in  the  case  of last-passage  percolation  described  in  the
introduction, we shall  call such a path a  \emph{geodesic}.  Let ${\cal N}_n$
denote the number  of collected points along such an  optimal path (which need
not be unique).  It is a well known fact (see \cite{AD}) that

\begin{equation}
  \label{2.1}
  \lim_{n\to + \infty} \frac{E^{(2)}{\cal N}_n}{n} = 2.
\end{equation}

In a remarkable paper \cite{J2}, K.   Johansson showed that for this model the
transversal fluctuations  of the geodesics  are of order $n^{2/3}$.  A closely
related problem is considered in \cite{BDJ}.

We now  modify the original  model in the  following way: Assume  in addition,
that  on the  main  diagonal  $y=x$ there  is  an independent  one-dimensional
Poisson point process  of intensity $\lambda^{(1)}=\lambda$, and let  $P^{(1)}$ be its
distribution.  We  will denote by $\Omega^{(1)}$  (resp. $\Omega^{(1)}_n$) the
configuration  space of  this process  (resp.\  of its  restriction to  ${\cal
  Q}_n$).   Finally  if  $\omega^{(1)}$  is  a realization  of  $P^{(1)}$  and
$\omega^{(2)}$ one of $P^{(2)}$, let $\omega = \omega^{(2)} \cup \omega^{(1)}$
be their union, and $\omega_n$ the restriction of $\omega$ to ${\cal Q}_n$.

What can  now be said about an  optimal directed path starting  at $(0,0)$ and
ending at $(n,n)$? If $\lambda\gg 1$,  clearly the geodesic will stay close to
the diagonal and
\begin{equation}\label{2.2}
  \lim_{n\to\infty}\frac{E^{(1)}\times          E^{(2)}          {\cal
  N}_n}{n}=:e(\lambda)>2.
\end{equation}
Instead of (\ref{2.2}) in our context it will be more convenient to use a more
geometric notion.  We say that the directed polymer is pinned (with respect to
the diagonal)  if $P^{(1)}\times P^{(2)}$-almost surely, the  number of visits
made by geodesics to $\{(u,u)\colon 0\le u\le t\}$ is of order $\gamma t$, for
some $\gamma>0$, for all $t$  large enough.  One expects the ``energy'' notion
(\ref{2.2}) and geometric  notion of pinning to be identical,  but this is yet
another point which remains to be proved.

$e(\lambda)$ is non-decreasing and there is a critical value $\lambda_\mathrm{c}$,
where   $e$  hits   the  value   $2$.   The   same  arguments   which  predict
$\lambda_\mathrm{c}=1$ in  the case of the  TASEP yield $\lambda_\mathrm{c}=0$
in the  case of our model.  Thus  any extra Poisson points  along the diagonal
are  expected to  pin  the directed  polymer.   Such a  behavior is  extremely
delicate, and the  answer depends on the nature and  behavior of the geodesics
in the initial, unperturbed system.  Very little is known, even on a heuristic
level, when  the underlying measure governing  the behavior of  the polymer is
not ``nice'' (with  a kind of Markov property, for  example a simple symmetric
random walk).  Our criteria presented  below give partial but rigorous answers
as to whether $\lambda_{\mathrm c}$ is strictly positive or not.

In passing let us note that for a symmetric environment pinning can be proved,
at  least on  the level  of  $e(\lambda)$ \cite{ImSa}.   Symmetric means  that
$P^{(2)}$ is concentrated on  point configurations which are symmetric relative
to   the   diagonal.   In   this   case  $\lambda_\mathrm{c}=1$,   \emph{i.e.}
$e(\lambda)=2$  for $0\leq\lambda\leq 1$  and $e(\lambda)>2$  for $\lambda>1$.
Indeed in the symmetric case, the  system is amenable to exact computations in
terms of Fredholm determinants; a trace of the simplification can also be seen
in our later discussion (see Section~\ref{pinning}).

\subsection{Construction of the broken-lines process}
\label{broken}

Hammersley gave a representation of the longest increasing subsequence problem
for a random  permutation in terms of broken lines built  from a Poisson point
process in the positive quadrant  (we describe the construction in some detail
below).  The length of the longest  increasing subsequence can then be seen to
be the number of lines which  separate the points $(0,0)$ and $(n,n)$ from one
another.  The  purpose of that  representation is to obtain  a superadditivity
property which easily implies the  existence of the limit in~(\ref{2.1}) --- but
doesn't specify its value.  It is  a very convenient formalism, which was used
in \cite{PS1,PS2,PS3} and \cite{SSV1,SSV2}.

The broken  line process $\Gamma_S$ in a  finite domain $S$ can  be defined as
the  space-time  trace   of  some  particle  system  with   birth,  death  and
immigration.  For convenience  we  rotate the  whole  picture by  an angle  of
$\pi/4$ clockwise,  so that the geodesic  is restricted to never  have a slope
which  is  larger  than  $1$  in  absolute value.  In  what  follows  we  will
consistently use the letters $t$ and $x$ for the first and second coordinates,
respectively, in  the rotated picture; we  will refer to $t$  as ``time'' (the
reason for that will become clear shortly).  The geodesics can then be seen as
curves of space-type (using the usual language of general relativity).

\begin{figure}
  \centering
  \begin{picture}(0,0)%
\includegraphics{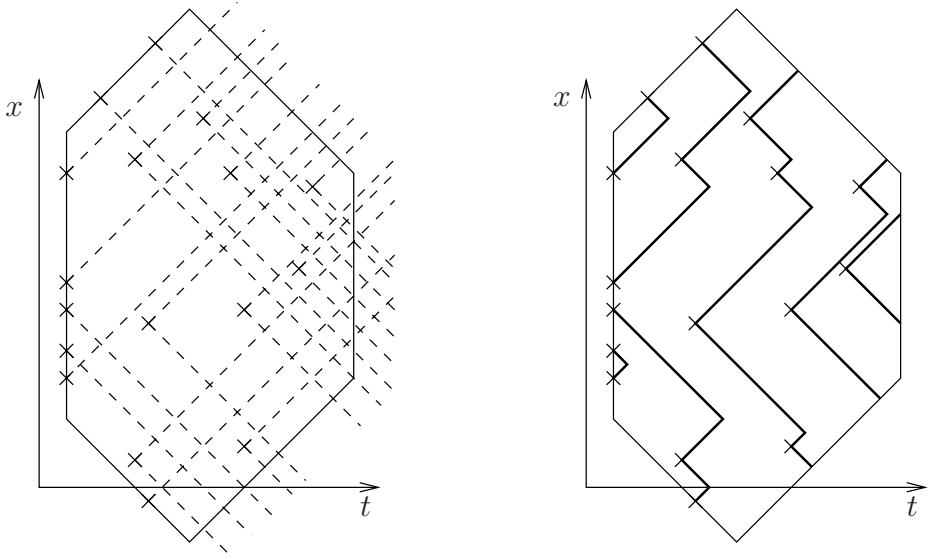}%
\end{picture}%
\setlength{\unitlength}{4144sp}%
\begingroup\makeatletter\ifx\SetFigFont\undefined%
\gdef\SetFigFont#1#2#3#4#5{%
  \reset@font\fontsize{#1}{#2pt}%
  \fontfamily{#3}\fontseries{#4}\fontshape{#5}%
  \selectfont}%
\fi\endgroup%
\begin{picture}(5416,3337)(361,-2441)
\put(361,209){\makebox(0,0)[rb]{\smash{\SetFigFont{12}{14.4}{\rmdefault}{\mddefault}{\updefault}{$x$}%
}}}
\put(2411,-2184){\makebox(0,0)[b]{\smash{\SetFigFont{12}{14.4}{\rmdefault}{\mddefault}{\updefault}{$t$}%
}}}
\put(5683,-2184){\makebox(0,0)[b]{\smash{\SetFigFont{12}{14.4}{\rmdefault}{\mddefault}{\updefault}{$t$}%
}}}
\put(3638,270){\makebox(0,0)[rb]{\smash{\SetFigFont{12}{14.4}{\rmdefault}{\mddefault}{\updefault}{$x$}%
}}}
\end{picture}
  \caption{The construction of a Hammersley process.}
  \label{fig:hammersley}
\end{figure}

Let   $S$    be   the   planar,    bounded   domain   defined    (cf.
Figure~\ref{fig:hammersley}) as
\begin{equation}\label{e2.3}
  S := \{(t,x): \, t_0 < t < t_1, \; g_{-}(t) < x < g_{+}(t) \},
\end{equation}
where $t_0 < t_1$ are given points  and $g_{-}(t) < g_{+}(t)$, $t_0 < t < t_1$
are  piecewise linear  continuous  functions such  that,  for some  $t_{01}^+$
(resp.  $t_{01}^-$)  in  the  interval $[t_0,t_1]$,  $g_{+}(t)$  (respectively
$g_{-}(t)$) increases (respectively,  decreases) on $(t_0,\, t_{01}^+)$ (resp.
$(t_0,\,  t_{01}^-)$) and decreases  (respectively, increases)  on $(t_{01}^+,\,
t_{1})$ (resp.  $(t_{01}^-,\, t_{1})$), always  forming an angle of  $\pm \pi/4$
with the $t$-axis.

Consider   four  independent   Poisson  processes   $\Pi_{0,+}$,  $\Pi_{0,-}$,
$\Pi_{+}$ and $\Pi_{-}$  on the boundary of $S$.   The processes $\Pi_{0,\pm}$
are supported on the  leftmost vertical boundary component $\bigtriangleup_0 S
:= \{ (t_0  ,x): \, g_{0,-} < x <  g_{0, +} \}$ of $S$, where  $ g_{0,\pm } :=
g_{\pm}(t_0)$,  and they  both have  intensity  $\sqrt{\lambda^{(2)}/2}$.  The
process $\Pi_{+}$ is defined on the ``northwest'' boundary $\bigtriangleup_+ S
:= \{ (t ,g_{0, +}(t)): \, t_0 <  t < t_{01}^+ \}$, and the process $\Pi_{-}$ is
defined  on the  ``southwest'' boundary  $\bigtriangleup_- S  := \{  (t ,g_{0,
  -}(t)): \,  t_0 <  t < t_{01}^-  \}$; their  intensities (with respect  to the
length element  of $\bigtriangleup_{\pm} S$)  are both $\sqrt{\lambda^{(2)}}$.
Finally,  let  $\Pi_{\text{in}}$ be  a  Poisson  point  process in  $S$,  with
intensity $\lambda^{(2)}$, and independent of the previous four.

Following the general  strategy for the definition of  Markov polygonal fields
of \cite{AS}, we  define a broken line process as follows.   Each point of the
Poisson process $\Pi_{\text{in}}$ is the point of birth of two particles which
start moving in opposite directions,  \emph{i.e.} with velocities $+1, \, -1$.
At  each  random point  of  $\Pi_{0,+}$, $\Pi_-$  a  particle  is born  having
velocity  $+1$.  Similarly  at each  random point  of $\Pi_{0,-}$,  $\Pi_+$, a
particle  is born  having  velocity  $-1$. All  particles  move with  constant
velocity until two  of them collide, after which  both colliding particles are
annihilated (see Figure~\ref{fig:hammersley}).

The  state space ${\cal  X}_S$ of  the process  $\Gamma_S$ is  the set  of all
finite collections $(\gamma_1, \dots , \gamma_k)$ (including the empty one) of
disjoint ``broken lines''  $\gamma_j$ inside $S$.  By a broken  line in $S$ we
mean the graph $\gamma = \{ (t,x) \in S: \, t = \widetilde \gamma(x)\}$ in $S$
of a continuous and piecewise linear function $\widetilde \gamma$, with slopes
all $\pm1$.   Let ${\cal P}_S$  denote the probability distribution  on ${\cal
  X}_S$ corresponding to the broken line process $\Gamma_S$ defined above.

Let $S^\prime  \subset S^{\prime \prime}$ be  two bounded domains  of the form
of~(\ref{e2.3}), and let ${\cal P}_{S^\prime}$ and ${\cal P}_{S^{\prime \prime}}$
be the  probability distributions of  the broken line processes  in $S^\prime$
and  $S^{\prime  \prime}$,   respectively.   The  probability  measure  ${\cal
  P}_{S^{\prime  \prime}}$   on  ${\cal  X}_{S^{\prime   \prime}}$  induces  a
probability  measure ${\cal  P}_{S^{\prime  \prime}}|_{{S^\prime}}$ on  ${\cal
  X}_{S^  \prime}$,  which  is  the  distribution of  the  restricted  process
$\Gamma_{S^{\prime \prime}}\cap  S^{\prime}$. Then, by  the choice we  made of
the  boundary  conditions,  the  following  consistency  property  holds  (see
\cite{AS}):
\begin{equation}
  {\cal P}_{S^{\prime \prime}}|_{{S^\prime}} = {\cal
    P}_{{S^\prime}}. \label{cdg}
\end{equation}
This  guarantees the  existence  of  the broken  line  process $\Gamma  \equiv
\Gamma_{\mathbb R^2}$ on $\mathbb R^2$,  its distribution ${\cal P}$ on ${\cal
  X}_{\mathbb R^2}$ being  such that for every $S$ of  the above shape, ${\cal
  P}|_{{S}} = {\cal P}_{{S}}$.  Moreover  ${\cal P}$ is invariant with respect
to the translations of $\mathbb R^2$.

\begin{rem}

The  same description  holds for  the  polynuclear growth  (PNG) model,  which
describes  a crystal  growing layer  by layer  on a  one-dimensional substrate
through  the  random  deposition  of  particles.  They  nucleate  on  existing
plateaus  of the  crystal,  forming  new islands.   In  an idealization  these
islands spread laterally with constant speed by steady condensation of further
material  at the  edges of  the islands.   Adjacent island  of the  same level
coalesce upon meeting and on the top of the new levels further islands emerge.

Observe that a path $\pi$ which  can move only in the northeast-southeast cone
(\emph{i.e.}  a  path of space type)  can collect at most  one initial Poisson
point from each  broken line.  In other words,  the broken lines ``factorize''
the points of the configuration $\omega$, in such a way that it tells us which
points cannot  be collected by the same  path, and that the  maximal number of
points  is  bounded  by the  number  of  broken  lines  which lie  in  between
start-point and end-point of the path $\pi$.

In fact,  the lines also  provide an explicit  construction of a  geodesic, as
follows:  Start at point  $(t_1,x)$ on  the right  boundary and  move leftward
until you  meet a line,  then follow  this line until  you arrive at  a point,
which you can  collect.  Then start moving leftward again  until you collect a
point in the  second-to-last broken line, and so on.   The number of collected
points  is then  essentially equal  to the  number of  broken lines,  though a
little care needs  to be taken as far as boundary  conditions are concerned if
this comparison is to be made completely formal. This observation led to a new
proof  of~(\ref{2.1})  in \cite{AS}.   It  is also  the  starting  point of  our
argument.
\end{rem}

\subsection{Essential and non-essential points}

We now return to the question asked  at the beginning of this section: How are
extra  added points  affecting the initial  system?  We  begin with  a  few purely
deterministic observations and statements.  We will need some extra notations:
Given any configuration $\widetilde  \omega_n$ (not necessarily sampled from a
Poisson process) of  points in ${\cal Q}_n $,  let $H(\widetilde \omega_n)$ be
the  number  of  broken  lines   produced  by  the  above  construction.   For
${\text{\bf x}}=(x,t)$ and $A=  \{{\text{\bf x}}_1,\dots , \, {\text{\bf x}}_k
\}$ we will denote by $\widetilde \omega_n \cup {\text{\bf x}}$ or $\widetilde
\omega_n \cup A$ the configuration  obtained from $\widetilde \omega_n$ by the
addition of the  points ${\text{\bf x}}_1$, \dots, ${\text{\bf  x}}_k$, and by
$\Gamma_n (\widetilde \omega_n)$ the associated configuration of broken lines.

\begin{proposition}[Abelian property, see \cite{SSV2}]
  \label{abelian1}
  For  any choice  of  $\widetilde \omega_n$,  ${\text{\bf x}}_1$  and
  ${\text{\bf x}}_2$ we have
  \begin{equation}
    \Gamma_n \left(  \widetilde \omega_n \cup  \left( {\text{\bf x}}_1
    \cup   {\text{\bf  x}}_2   \right)  \right)   =   \Gamma_n  \left(
    \left(\widetilde  \omega_n  \cup  {\text{\bf  x}}_1  \right)  \cup
    {\text{\bf  x}}_2  \right)   =  \Gamma_n  \left(  \left(\widetilde
      \omega_n  \cup {\text{\bf  x}}_2 \right)  \cup  {\text{\bf x}}_1
    \right).
  \end{equation}
\end{proposition}

\begin{proposition}[Monotonicity, see \cite{SSV2}]
  For any  choice of $\widetilde  \omega_n$ and ${\text{\bf  x}}_1$ we
  have
  \begin{equation}
    H (\widetilde  \omega_n )  \, \le \,  H (\widetilde  \omega_n \cup
    {\text{\bf x}}_1 )\, \le \, H (\widetilde \omega_n )+1.
  \end{equation}
\end{proposition}

\begin{definition}
  \label{essent}
  Given  $\widetilde  \omega_n$ and  ${\text{\bf  x}}  $  we say  that
  ${\text{\bf x}}  $ is \emph{essential for  $\widetilde \omega_n$} if
  $H(\widetilde  \omega_n  \cup {\text{\bf  x}}_1  )  = H  (\widetilde
  \omega_n )+1$.
\end{definition}

\begin{rem}
The above definition is domain-dependent: If $S_n \subset S_m$ are two domains
in  $\mathbb R^2$, and  $\widetilde \omega_n$,  $\widetilde \omega_m$  are the
restrictions  of  a  configuration  $\widetilde  \omega$ to  $S_n$  and  $S_m$
respectively, then  an extra  added point which  is essential  for $\widetilde
\omega_n$ might not be essential for $\widetilde \omega_m$.

If an added  point ${\text{\bf x}}$ is essential, its presence  is felt on the
boundary of the domain by the appearance of an extra broken line going outside
of  the area.   If  an added  point  ${\text{\bf x}}$  is  not essential,  its
presence can be  felt on the boundary or  not, but in any case  it will change
the local geometry of existing broken lines.

Speaking  informally,  the configuration  $\widetilde  \omega_n$ determines  a
partition of the domain into  two (possibly disconnected) regions $E$ and $B$,
such that any additional point chosen in $E$ will be essential for $\widetilde
\omega_n$, while  if it  is in $B$  it will  be not essential  for $\widetilde
\omega_n$.  It is  easy to  construct examples  of  configurations $\widetilde
\omega$  for which $E$  is empty  (\emph{i.e.}, that  are very  insensitive to
local changes).  On  the other hand $B$ is never empty  as soon as $\widetilde
\omega$  is not  empty. It  is also  easy to  give examples  of  the following
situations:
\begin{itemize}
\item  ${\text{\bf x}}_1  $ is  not  essential for  $\widetilde \omega_n$  and
  ${\text{\bf  x}}_2  $  is  not  essential  for  $\widetilde  \omega_n$,  but
  ${\text{\bf x}}_1  $ is essential  for $\widetilde \omega_n  \cup {\text{\bf
      x}}_2$  and ${\text{\bf x}}_2  $ is  essential for  $\widetilde \omega_n
  \cup {\text{\bf x}}_1$;
\item  ${\text{\bf  x}}_1  $   is  essential  for  $\widetilde  \omega_n$  and
  ${\text{\bf x}}_2 $ is  essential for $\widetilde \omega_n$, but ${\text{\bf
      x}}_2  $  is not  essential  for  $\widetilde  \omega_n \cup  {\text{\bf
      x}}_1$.
\end{itemize}
\end{rem}

\begin{rem}

In  order to  simplify our  explanations  and make  some concepts  as well
as the
notations   more  transparent  (and   lighter),  we   will  change   from  the
consideration  of   point-to-point  case  to   the  point-to-hyperplane  case,
\emph{i.e.} instead of  looking for a geodesic connecting  $(0,0)$ to $(t,0)$,
we will be  looking for an optimal path connecting $(t,0)$  to the line $x=0$.
This change is of purely ``pedagogical'' nature: All the ideas discussed above
and below are easily transferred  to the point-to-point case.  Nevertheless we
will not deny that it requires  some amount of additional work due to boundary
conditions.

The reader should  also not be surprised with our taking  of starting point as
$(t,0)$ and  moving backward to  the $x$-axis in the  point-to-hyperplane case
(see  the remark  at the  end of  Section~\ref{broken}).  Since  our broken
lines were constructed by drawing  the space-time trajectories of particles in
``forward  time'', the  information provided  by  the broken  lines about  the
underlying point configuration  is useful in the backward  direction, and thus
it  forces us  to construct  the geodesic  this way.  Conversely, in  order to
construct a forward geodesic, one could construct the broken lines on the same
point configuration but backward in time.

Due to that, the area to which  we will be restricting our process will be the
triangular area  $S_n$ enclosed by  segments connecting points  $(0,-n)$, $(0,
n)$ and  $(n,0)$.  If confusion doesn't  arise, we will  keep denoting related
quantities  by the  same sub-index  $n$ as  for the  square case,  for example
$\widetilde \omega_n$ will stay from now on for the configuration of points in
this triangular area.
\end{rem}

The next proposition is the crucial point in our construction.

\begin{proposition}[see \cite{SSV2}]
  \label{pin1}
  If  ${\text{\bf  x}} $  is  essential  for  $\widetilde \omega_n$  then  the
  point-to-plane   geodesic  in   configuration   $\widetilde  \omega_n   \cup
  {\text{\bf x}} $ has to collect point ${\text{\bf x}}$.
\end{proposition}

Again this  is a purely  deterministic statement, and  does not depend  on the
choice  of $\widetilde  \omega_n$ or  ${\text{\bf  x}}$.  Further  we will  be
considering only cases when extra points are added only along the $t$-axis.

\subsection{Propagation of influence}

Once an extra point ${\text{\bf x}}$ is added to the system, we need to update
the configuration  of broken lines. One  way to do  that is to redo  the whole
construction from scratch, \emph{i.e.} to  erase all the existing broken lines
and redraw  them using the algorithm  we described previously,  taking the new
point  into account.   It is  then natural  to ask  how much  the  new picture
differs  from  the old  one,  which  is not  perturbed  by  addition of  point
${\text{\bf x}}$.  It turns out that there is a very simple algorithm allowing
us to  trace all the  places of the  domain where the addition  of ${\text{\bf
    x}}$ will be felt, \emph{i.e.}  where local modification will be done.

Consider  an  augmented  configuration $\widetilde  \omega_n^\prime  =
\widetilde \omega_n \cup \{ {\text{\bf x}} \}$.  In order to see where
and how the broken lines  of $\Gamma_n (\widetilde \omega_n )$ will be
modified, we look at a  new interacting particle system, starting from
the points  of $\widetilde \omega_n^\prime$, but  with new interaction
rules:
\begin{enumerate}
\item  Particles  starting  from  the  points  of  $\widetilde  \omega_n$  are
  following the same rules as before, \emph{i.e.}  they are annihilated at the
  first time when  they collide with any other  particle. These particles will
  be called ``regular particles'';
\item The  two particles  which start from  the newly added  point ${\text{\bf
      x}}$  (also with  velocities $+1$  and $-1$)  will be  called ``superior
  particles'', and they obey different rules:
  \begin{enumerate}
  \item\label{annihila} Superior  particles annihilate if and  only if they  collide with each
    other;
  \item\label{annihil}  If  a superior  particle  collides  with some  regular
    particle, the  velocity of  the superior particle  changes to that  of the
    incoming  regular  particle,  which  is  annihilated  while  the  superior
    particle continues to move.
  \end{enumerate}
\end{enumerate}
We will  denote by $p_{{\text{\bf x}}}^+$ (resp.   $p_{{\text{\bf x}}}^-$) the
superior  particle which starts  from ${\text{\bf  x}}$ with  initial velocity
$+1$ (resp.  $-1$); The  space-time trajectories of $p_{{\text{\bf x}}}^+$ and
$p_{{\text{\bf  x}}}^-$  will  be  denoted  by  $\pi_{{\text{\bf  x}}}^+$  and
$\pi_{{\text{\bf x}}}^-$, respectively.

Observe  that if  any of  these two  trajectories leaves  the  triangular area
$S_n$, it  will never come back  to it. (Notice that  since superior particles
can change their  velocities during their evolution, both  particles can leave
the  triangular area  from the  same side.)   If $\pi_{{\text{\bf  x}}}^+$ and
$\pi_{{\text{\bf  x}}}^-$  intersect  inside  of  $S_n$,  then,  according  to
rule~\ref{annihila} the corresponding superior particles are annihilated.

The  path $\pi_{{\text{\bf  x}}}^+$ (resp.   $\pi_{{\text{\bf x}}}^-$)  can be
represented   as   an    alternating   sequence   of   concatenated   segments
$\pi_{{\text{\bf    x}}}^+    (1,+)$,    $\pi_{{\text{\bf    x}}}^+    (1,-)$,
$\pi_{{\text{\bf x}}}^+ (2,+)$, $\pi_{{\text{\bf  x}}}^+ (2,-), \dots $ (resp.
$\pi_{{\text{\bf    x}}}^-    (1,-)$,    $\pi_{{\text{\bf    x}}}^-    (1,+)$,
$\pi_{{\text{\bf x}}}^- (2,-)$, $\pi_{{\text{\bf x}}}^- (2,+), \dots $), where
each segment corresponds to the time interval between two consecutive velocity
changes of the superior particle $p_{{\text{\bf x}}}^+$ (resp., $p_{{\text{\bf
      x}}}^-$),  and  during which  its  velocity is  equal  to  $+1$ or  $-1$
according to the sign given as second argument in the notation.

\begin{figure}[htbp]
  \centering
 \begin{picture}(0,0)%
\includegraphics{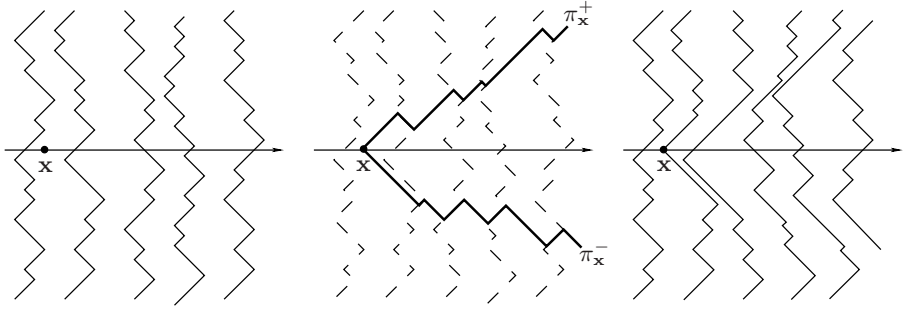}%
\end{picture}%
\setlength{\unitlength}{4144sp}%
\begingroup\makeatletter\ifx\SetFigFont\undefined%
\gdef\SetFigFont#1#2#3#4#5{%
  \reset@font\fontsize{#1}{#2pt}%
  \fontfamily{#3}\fontseries{#4}\fontshape{#5}%
  \selectfont}%
\fi\endgroup%
\begin{picture}(5397,1826)(214,-1694)
\put(433,-878){\makebox(0,0)[lb]{\smash{\SetFigFont{9}{10.8}{\familydefault}{\mddefault}{\updefault}{$\mathbf{x}$}%
}}}
\put(2336,-873){\makebox(0,0)[lb]{\smash{\SetFigFont{9}{10.8}{\familydefault}{\mddefault}{\updefault}{$\mathbf{x}$}%
}}}
\put(3575, 27){\makebox(0,0)[lb]{\smash{\SetFigFont{9}{10.8}{\rmdefault}{\mddefault}{\updefault}{$\pi_{\mathbf{x}}^+$}%
}}}
\put(3670,-1410){\makebox(0,0)[lb]{\smash{\SetFigFont{9}{10.8}{\rmdefault}{\mddefault}{\updefault}{$\pi_{\mathbf{x}}^-$}%
}}}
\put(4135,-873){\makebox(0,0)[lb]{\smash{\SetFigFont{9}{10.8}{\familydefault}{\mddefault}{\updefault}{$\mathbf{x}$}%
}}}
\end{picture}
  \caption{Propagation of influence.}
  \label{fig:propagation}
\end{figure}

The  trajectories of  $p_{{\text{\bf  x}}}^+$ and  $p_{{\text{\bf x}}}^-$  are
completely  determined by  $\widetilde \omega_n  \cup \{  {\text{\bf  x}} \}$.
Observe that each  time $p_{{\text{\bf x}}}^+$ changes its  velocity from $+1$
to  $-1$ (resp.   $p_{{\text{\bf x}}}^-$  changes  its velocity  from $-1$  to
$+1$), it  starts to move  along a segment  which also belongs to  some broken
line $\gamma_i$ from $\Gamma_n (\widetilde \omega_n)$, and when it changes its
velocity back  to $+1$, it leaves this  broken line, and moves  until the next
velocity flip, which happens exactly  when the superior particle collides with
the next broken line  $\gamma_{i+1}$ in $\Gamma_n (\widetilde \omega_n)$. This
gives  an  extremely  simple  rule  how  to  transform  $\Gamma_n  (\widetilde
\omega_n)$     into    $\Gamma_n     (\widetilde     \omega_n^\prime)$    (see
Figure~\ref{fig:propagation}):
\begin{itemize}
\item Erase all the $\pi_{{\text{\bf x}}}^\pm (j,\mp)$ (\emph{i.e.}, all parts
  of  the path of  the superior  particle which  are contained  in one  of the
  original broken  lines) to  obtain an intermediate  picture $\Gamma_n^\prime
  (\widetilde \omega_n)$;
\item  Add  all the  $\pi_{{\text{\bf  x}}}^\pm  (j,\pm)$  thus obtaining  the
  picture $\Gamma_n^{\prime \prime}(\widetilde \omega_n)$.
\end{itemize}

It is not hard to conclude that $\Gamma_n^{\prime \prime}(\widetilde \omega_n)
=   \Gamma_n  (\widetilde   \omega_n^\prime)$.  In   other  words   the  paths
$\pi_{{\text{\bf   x}}}^+$   and  $\pi_{{\text{\bf   x}}}^-$   show  how   the
``influence'' of  ${\text{\bf x}}$  spreads along the  configuration $\Gamma_n
(\widetilde \omega_n)$.  If both superior particles $p_{{\text{\bf x}}}^+$ and
$p_{{\text{\bf  x}}}^-$  collide inside  the  domain  $S_n$, the  trajectories
$\pi_{{\text{\bf x}}}^+$  and $\pi_{{\text{\bf x}}}^-$ close into  a loop, and
outside of  this loop the\vadjust{\goodbreak}  configuration $\Gamma_n (\widetilde  \omega_n)$ was
not modified,  \emph{i.e.} the  presence of ${\text{\bf  x}}$ was not  felt at
all.

\begin{definition}
  Given a  configuration $\omega$  of the underlying  Poisson process,  and an
  added point $\mathbf  x$ on the time axis, we  will denote by $\tau({\mathbf
    x};  \omega)$ (or  simply $\tau({\mathbf  x})$  if there  is no  confusion
  possible) the \emph{self-annihilation time} of the pair of particles created
  at $\mathbf x$, \emph{i.e.} the  time at which the paths $\pi_{\mathbf x}^-$
  and  $\pi_{\mathbf x}^+$  meet, if  such  a time  exists; let  $\tau(\mathbf
  x)=+\infty$ otherwise.
\end{definition}

In the  specific case we  are looking at,  $\tau(\mathbf x)$ is  almost surely
finite if $\lambda^{(2)}>0$, but it need  not be the case for other underlying
point processes.

\subsection{Interaction and Attractors}

Another  important step  in the  analysis of  the spread  of influence,  is to
understand how the influence paths interact with each other if we add multiple
points  ${\text{\bf  x}}_1, \dots  ,  {\text{\bf  x}}_\ell  $ to  the  initial
configuration. Proposition~\ref{abelian1} implies that  we can obtain the full
picture by  adding the points  one by one;  to simplify the notations,  in our
description of  the procedure we  will also use  the fact that  the additional
points will be placed along $t$-axis, though this is not essential.

Again, due to the presence of  time orientation, the nature of the interaction
between  influence paths  becomes  exposed in  a  more transparent  way if  we
proceed backward, \emph{i.e.}  if we  begin to observe the modifications first
when adding the  rightmost point, and then continue  progressively, adding the
points one by  one, moving leftward, each time checking  the effect created by
the newly added point.  For notational convenience let us index the new points
in the  backward direction,  \emph{i.e.}  ${\text{\bf x}}_i  = (t_i,  0)$ with
$t_1 > t_2 > \cdots$.

Applying the  construction described  in the previous  subsection successively
for each of the  new points, we obtain the following rules  for the updating a
configuration  with  multiple points  added:  Take  the initial  configuration
$\widetilde \omega_n$ and let $\widetilde \omega_n^{(r)} = \widetilde \omega_n
\cup \{ {\text{\bf x}}_{i}, 1\le i\le r \}$ be the modified configuration.  In
order  to  see  where  and  how  the broken  lines  of  $\Gamma_n  (\widetilde
\omega_n)$ will  be updated, consider  a new particle representation  built on
the  configuration  $\widetilde  \omega_n^{(r)}$,  and obeying  the  following
rules:
\begin{enumerate}
\item Regular  particles, starting from  the points of  $\widetilde \omega_n$,
  are annihilated as soon as they collide with any other particle;
\item The  particles starting  from the $({\text{\bf  x}}_i)_{1\le i  \le r}$,
  with velocities  $+1$ and $-1$  are denoted by $p_{{\text{\bf  x}}_i}^+$ and
  $p_{{\text{\bf  x}}_i}^-$ respectively;  again we  shall call  them superior
  particles. They behave as follows:
  \begin{enumerate}
  \item Whenever  two superior particle  of different types  collide (``$+-$''
    collision), they annihilate and both disappear;
  \item  If two  superior  particles of  the  same type  collide (``$++$''  or
    ``$-{}-$''  collision),  then  they  exchange  their  velocities  (elastic
    interaction) and continue to move;
  \item If a superior particle  collides with a regular particle, the velocity
    of the  superior particle changes  that of the incoming  regular particle;
    the regular particle is  annihilated, while the superior particle survives
    and continues to move.
  \end{enumerate}
\end{enumerate}

Denote  the  space-time  trajectories  of  superior  particles  $p_{{\text{\bf
      x}}_i}^+$ and $p_{{\text{\bf x}}_i}^-$ by $\pi_{{\text{\bf x}}_i}^+$ and
$\pi_{{\text{\bf  x}}_i}^-$,  respectively.  As  before,  each  pair of  paths
$\pi_{{\text{\bf x}}_i}^+$  and $\pi_{{\text{\bf x}}_i}^-$  can be represented
as an alternating sequence  of concatenated segments $\pi_{{\text{\bf x}}_i}^+
(1,+)$,  $\pi_{{\text{\bf x}}_i}^+  (1,-)$, $\pi_{{\text{\bf  x}}_i}^+ (2,+)$,
$\pi_{{\text{\bf x}}_i}^+ (2,-), \dots  $ or $\pi_{{\text{\bf x}}_i}^- (1,-)$,
$\pi_{{\text{\bf   x}}_i}^-    (1,+)$,   $\pi_{{\text{\bf   x}}_i}^-   (2,-)$,
$\pi_{{\text{\bf  x}}_i}^-  (2,+),  \dots   $,  respectively,  with  the  same
convention for the  sign of the velocities.  We are now  ready to complete the
set  of  rules  which  govern  the  transformation  of  $\Gamma_n  (\widetilde
\omega_n)$ in to $\Gamma_n (\widetilde \omega_n^{(r)})$:
\begin{itemize}
\item Erase  each segment $\pi_{{\text{\bf x}}_i}^\pm  (j,\mp)$ from $\Gamma_n
  (\widetilde \omega_n)$,  producing an intermediate  picture $\Gamma_n^\prime
  (\widetilde \omega_n)$;
\item   Add    the   segments   $\pi_{{\text{\bf    x}}_i}^\pm   (j,\pm)$   to
  $\Gamma_n^\prime  (\widetilde \omega_n)$,  thus  producing $\Gamma_n^{\prime
    \prime}(\widetilde \omega_n) = \Gamma_n (\widetilde \omega_n^{(r)})$.
\end{itemize}
(Here  as previously, two  $\pm$ in  the same  formula are  taken to  be equal
signs, while $\pm$ and $\mp$ in the same formula stand for opposite signs.)

Recall that  we are working  in the bounded  triangular domain $S_n$  with the
configuration $\widetilde  \omega_n \cup \{ {\text{\bf x}}_i,  1\le i\le r\}$,
where  $r$  is  the  number   of  added  points.   By  ${\text{\bf  f}}_i^+  =
(t_{{\text{\bf f}}_i^+},  x_{{\text{\bf f}}_i^+})$ and $  {\text{\bf f}}_i^- =
(t_{{\text{\bf  f}}_i^-},   x_{{\text{\bf  f}}_i^-})$  we   shall  denote  the
end-points   of   the   influence   paths   $\pi_{{\text{\bf   x}}_i}^+$   and
$\pi_{{\text{\bf  x}}_i}^-$ ---  they can  be points  where  the corresponding
paths  exit  the triangular  domain,  or  points  where a  ``$+-$''  collision
happens, in which case the two corresponding end-points are equal.

Besides,  let   ${\text{\bf  r}}_i^+   =  (t_{{\text{\bf  r}}_i^+},   0)$  and
${\text{\bf  r}}_i^- =  (t_{{\text{\bf r}}_i^-},  0)$ be  the points  of first
return  to   the  $t$-axis  of   the  paths  $\pi_{{\text{\bf   x}}_i}^+$  and
$\pi_{{\text{\bf x}}_i}^-$, respectively, and define
\begin{equation}
  \widehat  t_i  :=  \min  \{  t_{{\text{\bf  f}}_i^+},  t_{{\text{\bf
        f}}_i^-}, t_{{\text{\bf r}}_i^+}, t_{{\text{\bf r}}_i^-} \}.
\end{equation}
Last,  let ${\text{\bf  e}}_i^+ :=  \pi_{{\text{\bf x}}_i}^+  \cap \{(\widehat
t_i,  x), \, x\in  \mathbb R  \}$ and  ${\text{\bf e}}_i^-  := \pi_{{\text{\bf
      x}}_i}^- \cap \{(\widehat t_i, x), \,  x\in \mathbb R \}$, and denote by
$\widehat  \pi_{{\text{\bf x}}_i}^{+}$, $\widehat  \pi_{{\text{\bf x}}_i}^{-}$
the parts  of $\pi_{{\text{\bf  x}}_i}^{+}$ and $  \pi_{{\text{\bf x}}_i}^{-}$
lying between  ${\text{\bf x}}_i$  and ${\text{\bf e}}_i^+$,  and respectively
between ${\text{\bf x}}_i$ and ${\text{\bf e}}_i^-$.

\begin{figure}[htbp]
  \centering
  \begin{picture}(0,0)%
\includegraphics{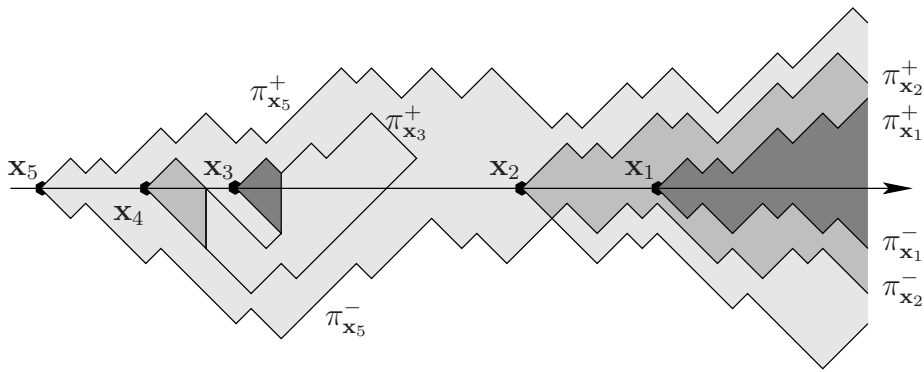}%
\end{picture}%
\setlength{\unitlength}{4144sp}%
\begingroup\makeatletter\ifx\SetFigFont\undefined%
\gdef\SetFigFont#1#2#3#4#5{%
  \reset@font\fontsize{#1}{#2pt}%
  \fontfamily{#3}\fontseries{#4}\fontshape{#5}%
  \selectfont}%
\fi\endgroup%
\begin{picture}(5424,2184)(6469,-4123)
\put(11701,-2401){\makebox(0,0)[lb]{\smash{\SetFigFont{12}{14.4}{\rmdefault}{\mddefault}{\updefault}{$\pi_{{\mathbf x}_2}^+$}%
}}}
\put(11701,-2671){\makebox(0,0)[lb]{\smash{\SetFigFont{12}{14.4}{\rmdefault}{\mddefault}{\updefault}{$\pi_{{\mathbf x}_1}^+$}%
}}}
\put(11701,-3391){\makebox(0,0)[lb]{\smash{\SetFigFont{12}{14.4}{\rmdefault}{\mddefault}{\updefault}{$\pi_{{\mathbf x}_1}^-$}%
}}}
\put(11701,-3661){\makebox(0,0)[lb]{\smash{\SetFigFont{12}{14.4}{\rmdefault}{\mddefault}{\updefault}{$\pi_{{\mathbf x}_2}^-$}%
}}}
\put(8371,-3841){\makebox(0,0)[lb]{\smash{\SetFigFont{12}{14.4}{\rmdefault}{\mddefault}{\updefault}{$\pi_{{\mathbf x}_5}^-$}%
}}}
\put(7921,-2491){\makebox(0,0)[lb]{\smash{\SetFigFont{12}{14.4}{\rmdefault}{\mddefault}{\updefault}{$\pi_{{\mathbf x}_5}^+$}%
}}}
\put(8731,-2671){\makebox(0,0)[lb]{\smash{\SetFigFont{12}{14.4}{\rmdefault}{\mddefault}{\updefault}{$\pi_{{\mathbf x}_3}^+$}%
}}}
\put(6661,-2941){\makebox(0,0)[rb]{\smash{\SetFigFont{12}{14.4}{\rmdefault}{\mddefault}{\updefault}{${\mathbf x}_5$}%
}}}
\put(7831,-2941){\makebox(0,0)[rb]{\smash{\SetFigFont{12}{14.4}{\rmdefault}{\mddefault}{\updefault}{${\mathbf x}_3$}%
}}}
\put(9541,-2941){\makebox(0,0)[rb]{\smash{\SetFigFont{12}{14.4}{\rmdefault}{\mddefault}{\updefault}{${\mathbf x}_2$}%
}}}
\put(10351,-2941){\makebox(0,0)[rb]{\smash{\SetFigFont{12}{14.4}{\rmdefault}{\mddefault}{\updefault}{${\mathbf x}_1$}%
}}}
\put(7291,-3211){\makebox(0,0)[rb]{\smash{\SetFigFont{12}{14.4}{\rmdefault}{\mddefault}{\updefault}{${\mathbf x}_4$}%
}}}
\end{picture}
  \caption{Attractors (shaded).}
  \label{fig:attractors}
\end{figure}

\begin{definition}
  \label{attr}
  Let  $J_i$ be  the (random)  Jordan  curve starting  at ${\text{\bf  x}}_i$,
  following  the  path $\widehat  \pi_{{\text{\bf  x}}_i}^-$  until the  point
  ${\text{\bf  e}}_i^-$,  then  the  vertical  line  $t=\widehat  t_i$  up  to
  ${\text{\bf e}}_i^+$, and then  the path $\widehat \pi_{{\text{\bf x}}_i}^+$
  backward until it  comes back to ${\text{\bf x}}_i$.   The domain bounded by
  $J_i$ will  be called the  \emph{attractor of the point  ${\text{\bf x}}_i$}
  and denoted by ${\cal  A}_i$ (see Figure~\ref{fig:attractors}).  The part of
  its  boundary  which is  contained  in  $\widehat \pi_{{\text{\bf  x}}_i}^+$
  (resp.  $\widehat \pi_{{\text{\bf x}}_i}^-$) will be called the \emph{upper}
  (resp.  \emph{lower}) \emph{boundary} of the attractor.
\end{definition}

It is important  to understand how attractors are affected  by one another, in
order to  give a convenient description  of the whole  augmented process.  The
key remark  is the  following: Informally speaking,  a superior particle  of a
given  type is  not affected  by an  older one  (\emph{i.e.}  one  with larger
index)  of the  same  type.  Indeed  this  is a  consequence  of the  previous
construction,    and    especially    of    the    Abelian    property    (see
Proposition~\ref{abelian1}).

Of course this does not mean  that the attractors ${\mathcal A}_j$ for $j<j_0$
do not  change when  ${\mathbf x}_{j_0}$ is  added, since it  remains possible
that  a ``$+-$''-collision  happens.  In  that  case, we  get some  kind of  a
monotonicity property, for  the statement of which some additional notation
will  be  needed.   Recall   that  $\tau(\mathbf  x;\omega)$  stands  for  the
self-annihilation  time  of  a  particle  \emph{in  the  underlying  scenery},
\emph{i.e.}  ignoring the effect of  the other new particles, both younger and
older.  We will denote by $\tau(\mathbf x, \mathbf y;\omega)$ the annihilation
time of the  set of superior particles born from $x$  and $y$, \emph{i.e.} the
last time at  which any of the corresponding four  superior particles is still
alive.  Such an annihilation can happen in one of four ways:
\begin{itemize}
\item Flat:  $\mathbf x$ is born,  the two particles  thus created annihilate,
  then $\mathbf y$ is born and its two particles collide;
\item  Embedded:  $\mathbf  x$  is  born, then  $\mathbf  y$  appears  between
  $p_{\mathbf  x}^-$ and  $p_{\mathbf x}^+$,  then the  particles  issued from
  $\mathbf y$ collide, then so do those issued from $\mathbf x$;
\item  Parallel: $\mathbf  x$ is  born, then  $\mathbf y$  appears  outside of
  $(p_{\mathbf x}^-,p_{\mathbf x}^+)$, then the particles issued from $\mathbf
  y$ collide, then so do those issued from $\mathbf x$;
\item  Crossed:  $\mathbf  x$ is  born,  then  $\mathbf  y$ (also  outside  of
  $(p_{\mathbf x}^-,p_{\mathbf x}^+)$), then one particle issued from $\mathbf
  x$ annihilates the particle of the  other type coming from $\mathbf y$, then
  the remaining two collide.
\end{itemize}

The combinatorics  become much  more involved when  more particles  are added;
nevertheless, it is possible (if a  bit technical if a formal proof is needed,
see~\cite{BS}) to show the following:
\begin{proposition}[Monotonicity of the influence]
  \label{vno}
  For any two added points $\mathbf  x$ and $\mathbf y$, we have the following
  inequality:
  $$\tau(\mathbf x,  \mathbf y;\omega) \ge  \mathrm{Max} (\tau(\mathbf
  x, \omega),  \tau(\mathbf y;\omega));$$ and more  generally the annihilation
  time of the union  of two finite families of added points  is at least equal
  to the larger of the two annihilation times of the parts.
\end{proposition}

In  the flat, embedded  and parallel  cases, the  monotonicity extends  to the
shapes of  the attractors  (the attractor  of $\mathbf y$  in the  presence of
$\mathbf  x$ contains the  one without);  there is  true reinforcement  in the
embedded case, in that the inclusion is strict as soon as there is a ``$++$''-
or   ``$-{}-$''-collision.   This   is  not   always  the   case   in  crossed
configurations (cf.~Fig.~\ref{fig:attractors}, where the attractor of $\mathbf x_3$
is  shortened  by the  addition  of  $\mathbf x_4$),  which  leads  us to  the
following definition:
\begin{definition}
  \label{atcon}
  We  say that  two  attractors  ${\cal A}_i$  and  ${\cal A}_j,  i  > j$  are
  \emph{connected} if there exists a sub-sequence  $j= i_0 < i_1 < i_2 < \dots
  < i_k =  i$ such that ${\text{\bf x}}_{i_r} \in  {\cal A}_{i_{r+1}}$ for all
  $0  \le r  <  k$. We  will  call $(i_j)_{0\le  j\le  k}$ a  \emph{connecting
    subsequence} between ${\cal A}_i$ and ${\cal A}_j$.
\end{definition}


Observe that  if $i>j>k$,  if ${\cal  A}_i$ is connected  to ${\cal  A}_j$ and
${\cal A}_j$ is  connected to ${\cal A}_k$, then ${\cal  A}_i$ is connected to
${\cal  A}_k$.   Nevertheless, due  to  the  presence  of orientation  in  the
temporal direction, the above implication  generally does not hold without the
condition $i>j>k$.

Our construction immediately implies the following:

\begin{proposition}
  \label{crab}
  If  ${\cal A}_i$  is  connected to  ${\cal  A}_j, \;  i  > j$,  then
  $\widehat t_i \ge \widehat t_j $.
\end{proposition}

\begin{corollary}
  \label{topdog}
  If  ${\cal A}_i$ is  connected to  ${\cal A}_j$,  $i>j$, and  the end-points
  ${\text{\bf  e}}_j^+$ and  ${\text{\bf e}}_j^-$  belong respectively  to the
  south-east and north-east boundaries of the triangular domain $S_n$, then so
  do ${\text{\bf e}}_i^+$ and ${\text{\bf e}}_i^-$.
\end{corollary}

\begin{corollary}
  \label{topdog1}
  Assume  that  ${\cal  A}_i$  is  connected  to  ${\cal  A}_j$,  $i>j$,  with
  connecting  subsequence  $(i_s)_{0\le s\le  k}$:  If  ${\text{\bf x}}_j$  is
  essential for the  configuration $\widetilde \omega_n$, then so  are all the
  ${\text{\bf x}}_{i_s}$, $1 \le s \le k$.
\end{corollary}

\subsection{Pinning of the geodesics}
\label{pinning}

We  now return  to  our original  problem.  Observe that  if,  for some  fixed
configuration $\widetilde  \omega_n$ in  $S_n$, we pick  a realization  of the
points $({\text{\bf x}}_i)_{1\le i\le v}$ in  such a way that ${\cal A}_1$ and
${\cal  A}_v$ are  connected (say),  with connecting  subsequence $(i_j)_{0\le
  j\le  k}$,  then  all  the  ${\text{\bf x}}_{i_j}$  must  be  essential  for
$\widetilde  \omega_n$, and  therefore  the point-to-plane  geodesics for  the
configuration $\widetilde \omega_n \cup \left\{ {\text{\bf x}}_i \right\}$ has
to visit all the ${\text{\bf x}}_{i_j}$.

In the new formalism, the original  question of whether, for any given density
$\lambda^{(1)}>0$ of  the one-dimensional Poisson point  process, the limiting
value in  (\ref{2.1}) is  increased, becomes equivalent  to the  following: Is
there a  positive $\delta$  such that,  with high probability  as $n$  goes to
infinity, at least a fraction $\delta$ of the newly added points are essential
for $\omega^{(2)}$?

This question is more complicated than  simply whether there exists a chain of
pairwise connected attractors spanning from  the left to the right boundary of
the domain: Indeed,  such a chain does not  necessarily have positive density.
In  the next  Section we  also mention  some of  the  interesting mathematical
questions that arise in the construction.

It is not an easy task to understand how the attractors behave.  The fact that
the structure  of the influence paths  $\widehat \pi_{{\text{\bf x}}_{i+1}}^+$
and $\widehat  \pi_{{\text{\bf x}}_{i+1}}^-$ depends  on $\Gamma_n (\widetilde
\omega_n \cup \{ {\text{\bf x}}_j, 1\le j\le  i \} )$, but not on the $\mathbf
x_j$  for  $j>i+1$,  reduces the  problem  to  checking  whether none  of  the
influence    paths   $\widehat   \pi_{{\text{\bf    x}}_{i+1}}^+$,   $\widehat
\pi_{{\text{\bf x}}_{i+1}}^-$ hits the  $t$-axis before ${\text{\bf x}}_i$, in
which   case  ${\cal   A}_{i+1}$   is  connected   to   ${\cal  A}_{i}$   (see
Figure~\ref{fig:attractors}).

For a single  point ${\text{\bf x}}$ added to  the initial configuration, each
influence path $\pi_{{\text{\bf x}}}^+$, $\pi_{{\text{\bf x}}}^-$ has the same
statistical properties as what is  known as a ``second-class particle'' in the
framework of exclusion processes.  Since  in the definition of an attractor an
important role is played by the (possible) return times of the influence paths
to the $t$-axis, several things must be settled:

\begin{enumerate}
\item \label{firstquestion} The first return  time to the $t$-axis of a single
  influence path. It is believed (but remains a challenging open problem) that
  in the case of a  one-dimensional exclusion process, a second-class particle
  behaves super-diffusively. Though some bounds are available, and we know the
  mean deviations of the second class particle \cite{PS4}, they do not provide
  good control on return times;
\item  The joint  behavior of  the influence  paths  $\widehat \pi_{{\text{\bf
        x}}}^+$  and $\widehat  \pi_{{\text{\bf x}}}^-$.   Generally, it  is a
  complicated question too, but for our  purposes we need to have such control
  only up  to the first times  when one of  $\widehat \pi_{{\text{\bf x}}}^+$,
  $\widehat  \pi_{{\text{\bf x}}}^-$ returns  to the  $t$-axis. Before  such a
  time, both paths stay apart from each other, and some good mixing properties
  of the system come into play;  so the question reduces to how efficiently we
  control point~\ref{firstquestion}.
\end{enumerate}

The fact that the influence lines of ``younger'' points (with smaller indices,
\emph{i.e.}   sitting  more  to  the  right) repeal  the  influence  lines  of
``older''  points, leads  to  the following  observation:  Once the  attractor
${\cal  A}_{i+m}$ of  an older  point  reaches the  younger point  ${\text{\bf
    x}}_i$, then it  cannot end before the attractor  ${\cal A}_{i}$ ends.  If
the attractor ${\cal  A}_{i}$ ends before reaching the  next point ${\text{\bf
    x}}_{i-1}$,  then ${\cal  A}_{i+m}$  can still  go  forward, and  possibly
itself reach ${\text{\bf x}}_{i-1}$.  Observe  that at the time ${\cal A}_{i}$
ends,  the  boundaries of  ${\cal  A}_{i+m}$  are  necessarily at  a  positive
distance from the $t$-axis.

That, together  with the fact  that the evolution of  $p_{{\text{\bf x}}_i}^+$
and  $p_{{\text{\bf x}}_i}^-$  in the  slab $(t_i,  t_{i-1})\times  \mathbb R$
depends only  on $\Gamma_n (\widetilde  \omega_n)$ brings some notion  of week
dependence to the  system from one side, and the idea  of a ``re-start point''
from another.  This  reduces the study of percolation of  attractors to a more
general problem of one-dimensional, long-range, dependent percolation which we
formulate in the next section.  There  we also mention some of the interesting
mathematical questions that have arisen.

\section{Stick percolation}

In  this section  we introduce  two ``stick  percolation'' models,  which will
serve  as toy  models in  the study  of the  propagation of  influence  in the
broken-line model. In spite of  their apparent simplicity, these models can be
very useful studying effects of columnar defects and establishing (bounds for)
critical   values  for   asymptotic  shape   changes  for   some   well  known
one-dimensional growth systems (see \cite{BS}).

\subsection{Model 1: overlapping sticks}

Let $(x_i)_{i\in\mathbb N}$ be a Poisson point process of intensity $\lambda >
0$ on  the positive real line. We  call the points of  this process ``seeds'',
and assume that they are ordered, $x_0$ being the point closest to the origin.
To  each  seed  $x_i$  we  associate  a  positive  random  variable  $S_i$  (a
``length'')  and assume  that  the $(S_i)_{i\in\mathbb  N}$  are i.i.d.\  with
common distribution function $F$.

The  system we  consider is  then the  following: For  every  $i\in\mathbb N$,
construct the segment  $\widehat S_i=[x_i, x_i +S_i]$ (which  we will call the
$i$-th ``stick'').  We say that  the sticks $\widehat S_i$ and $\widehat S_j$,
$i<j$ are  \emph{connected} if  $x_j <  x_i + S_i$,  \emph{i.e.} if  they have
non-empty intersection; the  system \emph{percolates} if and only  if there is
an  infinite  chain  of  distinct,  pairwise  connected  sticks,  which  (with
probability $1$) is equivalent to saying that the union of the sticks contains
a half-line.

It is easy to  see that the system percolates with probability  $0$ or $1$ (it
is a tail event  for the obvious filtration); and in fact  there is a complete
characterization of both cases:

\begin{proposition}[\cite{BS}]
  \label{percolation1}
  Let  $R(x)=P(S_1>x)$  be the  tail of  the
  stick length distribution, and let $\varphi (x) = \int_0^x R(u)du$. Then the
  system percolates with probability $1$ if and only if
  \begin{equation}
    \int_0^{+\infty} e^{-\lambda \varphi (x)} dx < + \infty.
  \end{equation}
\end{proposition}

This leads us to the following definition:

\begin{definition}\label{stabil}
  The  distribution   $F$  governing   the  stick  process   is  said   to  be
  \emph{cluster-stable} if  the system percolates for every  positive value of
  $\lambda$.
\end{definition}

\begin{example}[Return times of random walks]\label{ex1}
One natural distribution for the length of the sticks, in view of the previous
construction, is the following: At time  $x_i$, start two random walks with no
drift (the specifics, \emph{e.g.} whether they are discrete or continuous time
walks, will not  matter at this point  --- for that matter we  could also take
two  Brownian motions),  one from  $+1$ and  the other  from $-1$.   Then, let
$x_i+S_i$ be the first time when these two walks meet.  It is well known that,
up  to multiplicative constants,  $P(S_i>t)$ behaves  as $t^{-1/2}$  for large
$t$, as  soon as the  walks are irreducible  and their step  distribution have
finite  variance.  It  is  easy to  check  that the  obtained distribution  is
cluster-stable.
\end{example}

\begin{example}\label{ex2}
Change the  previous example  a little,  as follows: For  every value  of $i$,
start a two-dimensional Brownian motion (or random walk) starting at ($x_{i,1}$),
and let $x_i+S_i$ be
the first  hitting point of  the axis by  this Brownian motion; but  erase the
stick  if $S_i<0$.   Then  the distribution  $F$  is Cauchy  restricted to  be
positive, and  its tail is  equivalent to $P(S_i>t)  \sim c/t$ as $t$  goes to
infinity.  In that  case, the stick length distribution  is not cluster-stable
for model  1, and there exists a  critical value $\lambda_c =  1/c$, such that
for $\lambda > \lambda_c$ the system  percolates while it does not if $\lambda
< \lambda_c$.
\end{example}
\subsection{Model 2: reinforced sticks}

There are several  ways to mimic the interaction  between ``funnels''.  One of
the most simple possibilities is to add  extra rules to Model 1 to account for
the interaction.  The  basic idea behind this modification, is  that if two or
more sticks from Model 1 overlap, then there is a certain reinforcement of the
system, which depends on how many sticks overlap, and then the whole connected
component is enlarged correspondingly.

One  way to  do that  can  be described  informally as  the following  dynamic
process.   First, see  each  stick $\widehat  S_i$  as the  flight  time of  a
particle $\pi_i$ born at time $x_i$. We want to model the fact that if $\pi_i$
wants  to land  when a  younger  one (say  $\pi_j$) is  still flying,  instead
$\pi_i$  ``bounces'' on  $\pi_j$;  $\pi_j$ on  the  other hand  should not  be
affected  by $\pi_i$,  if  the process  is  to look  like  the propagation  of
influence described in he previous section.

Assign to  each particle a  ``counter of chances''  $N_i$ which is set  to $1$
when the particle is born (formally it should be a function from $\mathbb R_+$
to $\mathbb  N$, and we let  $N_i(x_i)=1$).  Then, two things  can happen.  If
$x_i+S_i<x_{i+1}$, there is no interaction and $\pi_i$ dies at time $x_i+S_i$.
If on  the other hand  $x_i+S_i \ge x_{i+1}$,  at time $x_{i+1}$  the particle
$\pi_i$ gets a  ``bonus'', so that $N_i(x_{i+1})=2$; and  similarly, it gets a
bonus each time it passes above a seed point, so that $N_i(x_i+S_i-)-1$ is the
number of seed points in $\widehat  S_i$. Now when $\pi$ lands, its counter is
decreased by $1$, but if it is still positive the particle bounces on the axis
and  restarts using  an  independent copy  of  $S$.  If  all  the chances  are
exhausted before the closest seed is reached, $\pi_i$ gets killed.

Again,  we  may   ask  a  similar  question:  Given   $\lambda>0$,  for  which
distribution functions $F$ is the  probability for a given particle to survive
up to infinity positive?  It is obvious that even with $E(S) < +\infty$ we can
obtain infinite trajectories for certain  (large) values of $\lambda$. We will
call $F$ \emph{cluster-stable for the  reinforced process} if this happens for
every positive $\lambda$.

\begin{examp}
The Cauchy  distribution used in  the previous example cluster-stable  for the
reinforced model.

It is easy  to see (\emph{e.g.}  by a coupling argument)  that if the original
system percolates, the reinforced version (for the same value of $\lambda$ and
the   same  length  distribution   $F$)  percolates   too.  In   particular  a
cluster-stable  distribution  is cluster-stable  for  the reinforced  problem.
Nevertheless it is  an interesting open problem to  give full characterization
of distributions which are cluster-stable for the reinforced models.
\end{examp}

\paragraph{Acknowledgments.}

The authors wish  to thank K.  Alexander, H. Kesten and  D. Surgailis for many
hours  of fruitful  and clarifying  discussions scattered  over the  past four
years.  We also thank CBPF, IMPA, TU-M\"{u}nchen for hospitality and financial
support.


\begin{thebibliography}{99}

\bibitem{AD}
{\sc Aldous, D.} {\sc and} {\sc Diaconis, P.} (1995).
\newblock Hammersley's interacting particle process and longest increasing
  subsequences.
\newblock {\em Probab. Theory Related Fields\/}~\textbf{103},~2, 199--213.
\MR{1355056}

\bibitem{AS}
{\sc Arak, T.} {\sc and} {\sc Surgailis, D.} (1989).
\newblock Markov fields with polygonal realizations.
\newblock {\em Probab. Theory Related Fields\/}~\textbf{80},~4, 543--579.
\MR{980687}

\bibitem{BDJ}
{\sc Baik, J.}, {\sc Deift, P.}, {\sc and} {\sc Johansson, K.}
(1999).
\newblock On the distribution of the length of the longest increasing
  subsequence of random permutations.
\newblock {\em J.~Amer. Math. Soc.\/}~\textbf{12},~4, 1119--1178.
\MR{1682248}

\bibitem{BS}
{\sc Beffara, V.} {\sc and} {\sc Sidoravicius, V.} (2005)
\newblock {Effect of columnar defect on asymptotic shape of some growth
  processes.}
\newblock Preprint, in preparation.


\bibitem{CR}
{\sc Covert, P.} {\sc and} {\sc Rezakhanlou, F.} (1997).
\newblock Hydrodynamic limit for particle systems with non-constant speed
  parameter.
\newblock {\em J. Statist. Phys.\/}~\textbf{88},~1-2, 383--426.
\MR{1468390}

\bibitem{HTN}
{\sc Ha, M.}, {\sc Timonen, J.}, {\sc and} {\sc den Nijs, M.} (2003)
\newblock Queuing transitions in the asymmetric simple exclusion process.
\newblock {\em Phys. Rev. E}, 68:056122.


\bibitem{ImSa}
{\sc Sasamoto, T.} {\sc and} {\sc Imamura, T.} (2004).
\newblock Fluctuations of the one-dimensional polynuclear growth model in
  half-space.
\newblock {\em J. Statist. Phys.\/}~\textbf{115},~3-4, 749--803.
\MR{2054161}

\bibitem{JL2}
{\sc Janowsky, S.~A.} {\sc and} {\sc Lebowitz, J.~L.} (1994).
\newblock Exact results for the asymmetric simple exclusion process with a
  blockage.
\newblock {\em J. Statist. Phys.\/}~\textbf{77},~1-2, 35--51.
\MR{1300527}

\bibitem{J2}
{\sc Johansson, K.} (2000).
\newblock Transversal fluctuations for increasing subsequences on the plane.
\newblock {\em Probab. Theory Related Fields\/}~\textbf{116},~4, 445--456.
\MR{1757595}

\bibitem{L}
{\sc Liggett, T.~M.} (1999).
\newblock {\em Stochastic Interacting Systems: Contact, Voter and Exclusion
  Processes}. Grundlehren der Mathematischen Wissenschaften [Fundamental
  Principles of Mathematical Sciences], Vol. \textbf{324}.
\newblock Springer-Verlag, Berlin.
\MR{1717346}

\bibitem{PS1}
{\sc Pr{\"a}hofer, M.} {\sc and} {\sc Spohn, H.} (2000).
\newblock Statistical self-similarity of one-dimensional growth processes.
\newblock {\em Phys. A\/}~\textbf{279},~1-4, 342--352.
\MR{1797145}

\bibitem{PS4}
{\sc Pr{\"a}hofer, M.} {\sc and} {\sc Spohn, H.} (2002).
\newblock Current fluctuations for the totally asymmetric simple exclusion
  process.
\newblock In {\em In and Out of Equilibrium (Mambucaba, 2000)}. Progr. Probab.,
  Vol.~\textbf{51}. Birkh\"auser Boston, Boston, MA, 185--204.
\MR{1901953}

\bibitem{PS2}
{\sc Pr{\"a}hofer, M.} {\sc and} {\sc Spohn, H.} (2002).
\newblock Scale invariance of the {PNG} droplet and the {A}iry process.
\newblock {\em J. Statist. Phys.\/}~\textbf{108},~5-6, 1071--1106.
\MR{1933446}

\bibitem{PS3}
{\sc Pr{\"a}hofer, M.} {\sc and} {\sc Spohn, H.} (2004).
\newblock Exact scaling functions for one-dimensional stationary {KPZ} growth.
\newblock {\em J. Statist. Phys.\/}~\textbf{115},~1-2, 255--279.
\MR{2070096}

\bibitem{S}
{\sc Sepp{\"a}l{\"a}inen, T.} (2001).
\newblock Hydrodynamic profiles for the totally asymmetric exclusion process
  with a slow bond.
\newblock {\em J. Statist. Phys.\/}~\textbf{102},~1--2, 69--96.
\MR{1819699}

\bibitem{SSV1}
{\sc Sidoravicius, V.}, {\sc Vares, M.~E.}, {\sc and} {\sc
Surgailis, D.}
  (1999).
\newblock Poisson broken lines process and its application to {B}ernoulli first
  passage percolation.
\newblock {\em Acta Appl. Math.\/}~\textbf{58},~1--3, 311--325.
\MR{1734758}

\bibitem{SSV2}
{\sc Sidoravicius, V.}, {\sc Vares, M.~E.}, {\sc and} {\sc
Surgailis, D.} (2005).
\newblock Discrete broken line process and applications to the first passage
  percolation models.
\newblock  Preprint, \newblock in preparation.


\end{thebibliography}
\end{document}